\documentclass[12pt,shift]{article}

\newtheorem{thm}{Theorem}[section]
\newtheorem{cor}[thm]{Corollary}
\newtheorem{lem}[thm]{Lemma}

\begin{document}
\title{ {\bf The high-level error bound for shifted surface spline interpolation}}

\author{ {\bf Lin-Tian Luh\thanks{This work was supported by NSC 93-2115-M126-004.}} \\
Department of Mathematics, Providence University, \\ Shalu, Taichung, Taiwan \\ email:ltluh@pu.edu.tw \\phone:(04)26328001 ext. 15126 \\fax:(04)26324653}

\maketitle

\bigskip
Radial function interpolation of scattered data is a frequently used method for multivariate data fitting. One of the most frequently used radial functions is called shifted surface spline, introduced by Dyn, Levin and Rippa in \cite{Dy1} for $R^{2}$. Then it's extended to $R^{n}$ for $n\geq 1$. Many articles have studied its properties, as can be seen in \cite{Bu,Du,Dy2,Po,Ri,Yo1,Yo2,Yo3,Yo4}. When dealing with this function, the most commonly used error bounds are the one raised by Wu and Schaback in \cite{WS}, and the one raised by Madych and Nelson in \cite{MN2}. Both are $O(d^{l})$ as $d\rightarrow 0$, where $l$ is a positive integer and $d$ is the fill-distance. In this paper we present an improved error bound which is $O(\omega^{1/d})$ as $d\rightarrow 0$, where $0<\omega <1$ is a constant which can be accurately calculated.\\
{\bf Keywords}:radial basis function, shifted surface spline, error bound.\\
{\bf AMS subject classification}: 41A05, 41A15, 41A25, 41A30, 41A63, 65D10. 

\pagenumbering{arabic}

\setcounter{section}{0}

\section{Introduction}

Let $h$ be a continuous function on $R^{n}$ which is conditionally positive definite of order $m$. Given data $(x_{j}, f_{j}),\ j=1,\ldots , N$, where $X=\{ x_{1},\ldots , x_{N} \} $ is a subset of points in $R^{n}$ and the $f_{j}'s$ are real or complex numbers, the so-called $h$ spline interpolant of these data is the function $s$ defined by
\begin{equation}
  s(x)=p(x)+\sum_{j=1}^{N}c_{j}h(x-x_{j}),   \label{eq:11}
\end{equation}
where $p(x)$ is a polynomial in $P_{m-1}$ and $c_{j}'s$ are chosen so that
\begin{equation}
  \sum_{j=1}^{N}c_{j}q(x_{j})=0   \label{eq:12}
\end{equation}
for all polynomials $q$ in $P_{m-1}$ and
\begin{equation}
  p(x_{i})+\sum_{j=1}^{N}c_{j}h(x_{i}-x_{j})=f_{i},\ i=1,\ldots , N.   \label{eq:13}
\end{equation}
Here $P_{m-1}$ denotes the class of those polynomials of $R^{n}$ of degree $\leq m-1$.

It is well known that the system of equations (\ref{eq:12}) and (\ref{eq:13}) has a unique solution when $X$ is a determining set for $P_{m-1}$ and $h$ is strictly conditionally positive definite. For more details please see \cite{MN1}. Thus, in this case, the interpolant $s(x)$ is well defined.

We remind the reader that $X$ is said to be a determining set for $P_{m-1}$ if $p$ is in $P_{m-1}$ and $p$ vanishes on $X$ imply that $p$ is identically zero.

In this paper $h$ is defined by formula
\begin{eqnarray}
  h(x)&:=&(-1)^{m}(|x|^{2}+c^{2})^{\frac{\lambda}{2}}log(|x|^{2}+c^{2})^{\frac{1}{2}},\ \lambda \in Z_{+},\ m=1+\frac{\lambda}{2},\ c>0,     \nonumber \\
      &  & x\in R^{n},\ \lambda,\ n\ even,           
\end{eqnarray}
where $|x|$ is the Euclidean norm of $x$, and $\lambda,\ c$ are constants. In fact, the definition of shifted surface spline covers odd dimensions. For odd dimensions, it's of the form
\begin{eqnarray}
  h(x)&:=&(-1)^{\lceil \lambda -\frac{n}{2}\rceil }(|x|^{2}+c^{2})^{\lambda-\frac{n}{2}},\ n\ odd,\ \lambda \in Z_{+}=\{ 1,2,3,\ldots \}\  \nonumber \\
      &  &and\ \lambda > \frac{n}{2}.
\end{eqnarray}
However, this is just multiquadric whose exponential error estimates have already been constructed by Madych and Nelson in \cite{MN3}, and the calculation of the constant $\omega$ in $O(\omega^{\frac{1}{d}})$ can be found in \cite{Lu5}.Hence we will not discuss it. Instead, we will focus on even dimensions.
\subsection{A Bound for Multivariate Polynomials}

A key ingredient in the development of our estimates is the following lemma which gives a bound on the size of a polynomial on a cube in $R^{n}$ in terms of its values on a discrete subset which is scattered in a sufficiently uniform manner. We cite it directly from \cite{MN3} and omit its proof.
\begin{lem}
  For $n=1,2,\ldots ,$ define $\gamma_{n}$ by the formulae $\gamma_{1}=2$ and, if $n>1$, $\gamma_{n}=2n(1+\gamma_{n-1})$. Let $Q$ be a cube in $R^{n}$ that is subdivided into $q^{n}$ identical subcubes. Let $Y$ be a set of $q^{n}$ points obtained by selecting a point from each of those subcubes. If $q\geq \gamma_{n}(k+1)$, then for all $p$ in $P_{k}$ 
$$sup_{x\in Q}|p(x)|\leq e^{2n\gamma_{n}(k+1)}sup_{y\in Y}|p(y)|.$$
\end{lem}
\subsection{A Variational Framework for Interpolation}

The precise statement of our estimates concerning $h$ splines requires a certain amount of technical notation and terminology which is identical to that used in \cite{MN2}. For the convenience of the reader we recall several basic notions.

The space of complex-valued functions on $R^{n}$ that are compactly supported and infinitely differentiable is denoted by ${\cal D}$. The Fourier transform of a function $\phi$ in ${\cal D}$ is
$$\hat{\phi}(\xi)=\int e^{-i<x,\xi>}\phi(x)dx.$$
A continuous function $h$ is conditionally positive definite of order $m$ if $$\int h(x)\phi(x)\ast \tilde{\phi}(x)dx\geq 0$$ holds whenever $\phi=p(D)\psi $ with $\psi$ in ${\cal D}$ and $p(D)$ a linear homogeneous constant coefficient differential operator of order $m$. Here $\tilde{\phi}=\overline{\phi(-x)}$ and $\ast $ denotes the convolution product 
$$\phi_{1}\ast \phi_{2}(t)=\int \phi_{1}(x)\phi_{2}(t-x)dx.$$
As pointed out in \cite{MN2}, this definition of conditional positive definiteness is equivalent to that of \cite{MN1} which is generally used.

If $h$ is a continuous conditionally positive definite function of order $m$, the Fourier transform of $h$ uniquely determines a positive Borel measure $\mu$ on $R^{n}\backslash \{ 0\} $ and constants $a_{\gamma},\ |\gamma|=2m$ as follows:
 For all $\psi \in {\cal D}$
\begin{eqnarray}
  \int h(x)\psi(x)dx & = & \int \{ \hat{\psi}(\xi)-\hat{\chi }(\xi)\sum_{|\gamma|<2m}D^{\gamma}\hat{\psi}(0)\frac{\xi^{\gamma}}{\gamma !}\} d\mu (\xi)  \nonumber \\
                     &   & +\sum_{|\gamma|\leq 2m}D^{\gamma}\hat{\psi}(0)\frac{a_{\gamma}}{\gamma !},  
\end{eqnarray}
where for every choice of complex numbers $c_{\alpha},\ |\alpha|=m$,
$$\sum_{|\alpha|=m}\sum_{|\beta|=m}a_{\alpha+\beta}c_{\alpha}\overline{c_{\beta}}\geq 0.$$
Here $\chi$ is a function in ${\cal D}$ such that $1-\hat{\chi}(\xi)$ has a zero of order $2m+1$ at $\xi=0$; both of the integrals $\int_{0<|\xi|<1}|\xi|^{2m}d\mu(\xi),\ \int_{|\xi|\geq 1}d\mu(\xi)$ are finite. The choice of $\chi$ affects the value of the coefficients $a_{\gamma}$ for $|\gamma|<2m$.

Our variational framework for interpolation is supplied by a space we denote by ${\cal C}_{h,m}$. If 
$${\cal D}_{m}=\{ \phi\in {\cal D}: \int x^{\alpha}\phi(x)dx=0\ for\ all\ |\alpha|<m\} ,$$
then ${\cal C}_{h,m}$ is the class of those continuous functions $f$ which satisfy 
\begin{equation}
  \left| \int f(x)\phi(x)dx\right| \leq c(f)\left\{ \int h(x-y)\phi(x)\overline{\phi(y)}dxdy\right\} ^{\frac{1}{2}}
\end{equation}
for some constant $c(f)$ and all $\phi$ in ${\cal D}_{m}$. If $f\in {\cal C}_{h,m}$, let $\| f\| _{h}$ denote the smallest constant $c(f)$ for which (7) is true. Recall that $\| f\| $ is a semi-norm and ${\cal C}_{h,m}$ is a semi-Hilbert space; in the case $m=0$ it is a norm and a Hilbert space respectively.

\section{Main Results}

We first recall that the function $h$ defined in (4) is conditionally positive definite of order $m=1+\frac{\lambda}{2}$. This can be found in \cite{Dy2} and many relevant papers. Its Fourier transform \cite{GS} is 
\begin{equation}
  \hat{h}(\theta)=l(\lambda,n)|\theta|^{-\lambda-n}\tilde{{\cal K}}_{\frac{n+\lambda}{2}}(c|\theta|)
\end{equation} 
where $l(\lambda,n)>0$ is a constant depending on $\lambda$ and $n$, and $\tilde{{\cal K}}_{\nu}(t)=t^{\nu}{\cal K}_{\nu}(t)$, ${\cal K}_{\nu}(t)$ being the modified Bessel function of the second kind\cite{AS}. Then we have the following lemma.
\begin{lem}
  Let $h$ be as in (4) and $m$ be its order of conditional positive definiteness. There exists a positive constant $\rho$ such that 
\begin{equation}
  \int_{R^{n}}|\xi|^{k}d\mu(\xi)\leq l(\lambda,n)\cdot \sqrt{\frac{\pi}{2}}\cdot n\cdot \alpha_{n}\cdot c^{\lambda-k}\cdot \Delta_{0}\cdot \rho^{k}\cdot  k!
\end{equation}
for all integer $k\geq 2m+2$ where $\mu$ is defined in (6), $\alpha_{n}$ denotes the volume of the unit ball in $R^{n}$, $c$ is as in (4), and $\Delta_{0}$ is a positive constant.
\end{lem}
{\bf Proof}. We first transform the integral of the left-hand side of the inequality into a simpler form.
\begin{eqnarray}
  &      & \int_{R^{n}}|\xi|^{k}d\mu(\xi)  \nonumber \\
  & =    & \int_{R^{n}}|\xi|^{k}l(\lambda,n)\tilde{{\cal K}}_{\frac{n+\lambda}{2}}(c|\xi|)|\xi|^{-\lambda-n}d\xi \nonumber \\
  & =    & l(\lambda,n)c^{\frac{n+\lambda}{2}}\int_{R^{n}}|\xi|^{k-\frac{n+\lambda}{2}}\cdot {\cal K}_{\frac{n+\lambda}{2}}(c|\xi|)d\xi \nonumber \\
  & \sim & l(\lambda,n)c^{\frac{n+\lambda}{2}}\sqrt{\frac{\pi}{2}}\int_{R^{n}}|\xi |^{k-\frac{n+\lambda}{2}}\cdot \frac{1}{\sqrt{c|\xi|}\cdot e^{c|\xi|}}d\xi \nonumber \\
  & =    & l(\lambda,n)c^{\frac{n+\lambda}{2}}\cdot \sqrt{\frac{\pi}{2}}\cdot n\cdot \alpha_{n}\int_{0}^{\infty}r^{k-\frac{n+\lambda}{2}}\cdot \frac{r^{n-1}}{\sqrt{cr}\cdot e^{cr}}dr  \nonumber \\
  & =    & l(\lambda,n)c^{\frac{n+\lambda}{2}}\sqrt{\frac{\pi}{2}}\cdot n\cdot \alpha_{n}\cdot \frac{1}{\sqrt{c}}\int _{0}^{\infty}\frac{r^{k+\frac{n-\lambda -3}{2}}}{e^{cr}}dr  \nonumber \\
  & =    & l(\lambda,n)c^{\frac{n+\lambda}{2}}\sqrt{\frac{\pi}{2}}\cdot n\cdot \alpha_{n}\cdot \frac{1}{\sqrt{c}}\cdot \frac{1}{c^{k+\frac{n-\lambda -1}{2}}}\int_{0}^{\infty}\frac{r^{k+\frac{n-\lambda -3}{2}}}{e^{r}}dr  \nonumber \\
  & =    & l(\lambda,n)\sqrt{\frac{\pi}{2}}\cdot n\cdot \alpha_{n}\cdot c^{\lambda-k}\int _{0}^{\infty}\frac{r^{k'}}{e^{r}}dr \ where\ k'=k+\frac{n-\lambda -3}{2}. \nonumber 
\end{eqnarray}
Note that $k\geq 2m+2=4+\lambda$ implies $k'\geq \frac{n+\lambda+5}{2}>0$.

Now we divide the proof into three cases. Let $k''=\lceil k'\rceil $ which is the smallest integer greater than or equal to $k'$.

Case1. Assume $k''>k$. Let $k''=k+s$. Then 
$$\int_{0}^{\infty}\frac{r^{k'}}{e^{r}}dr\leq \int_{0}^{\infty}\frac{r^{k''}}{e^{r}}dr=k''!=(k+s)(k+s-1)\cdots (k+1)k!$$
and
$$\int_{0}^{\infty}\frac{r^{k'+1}}{e^{r}}dr\leq \int_{0}^{\infty}\frac{r^{k''+1}}{e^{r}}dr=(k''+1)!=(k+s+1)(k+s)\cdots (k+2)(k+1)k!.$$
Note that
$$\frac{(k+s+1)(k+s)\cdots (k+2)}{(k+s)(k+s-1)\cdots (k+1)}=\frac{k+s+1}{k+1}.$$
The condition $k\geq 2m+2$ implies that
$$\frac{k+s+1}{k+1}\leq \frac{2m+3+s}{2m+3}=1+\frac{s}{2m+3}.$$
Let $\rho=1+\frac{s}{2m+3}$. Then
$$\int_{0}^{\infty}\frac{r^{k''+1}}{e^{r}}dr\leq \Delta_{0}\cdot \rho^{k+1}\cdot (k+1)!$$
if $\int_{0}^{\infty}\frac{r^{k''}}{e^{r}}dr\leq \Delta_{0}\cdot \rho^{k}\cdot k!$. The smallest $k''$ is $k_{0}''=2m+2+s$ when $k=2m+2$. Now,
\begin{eqnarray}
  \int_{0}^{\infty}\frac{r^{k_{0}''}}{e^{r}}dr & = & k_{0}''!=(2m+2+s)(2m+1+s)\cdots(2m+3)(2m+2)! \nonumber \\
                                               & = & \frac{(2m+2+s)(2m+1+s)\cdots (2m+3)}{\rho^{2m+2}}\cdot \rho^{2m+2}\cdot (2m+2)! \nonumber \\
                                               & = & \Delta_{0}\cdot \rho^{2m+2}\cdot(2m+2)!   \nonumber \\ 
                                               &   &  where\ \Delta_{0}=\frac{(2m+2+s)(2m+1+s)\cdots (2m+3)}{\rho^{2m+2}}. \nonumber   
\end{eqnarray}
It follows that $\int_{0}^{\infty}\frac{r^{k'}}{e^{r}}dr\leq \Delta_{0}\cdot \rho^{k}\cdot k!$ for all $k\geq 2m+2$.

Case2. Assume $k''<k$. Let $k''=k-s$ where $s>0$. Then 
$$\int_{0}^{\infty}\frac{r^{k'}}{e^{r}}dr\leq \int_{0}^{\infty}\frac{r^{k''}}{e^{r}}dr=k''!=(k-s)!=\frac{1}{k(k-1)\cdots (k-s+1)}\cdot k!$$
and
\begin{eqnarray}
  \int_{0}^{\infty}\frac{r^{k'+1}}{e^{r}}dr & \leq & \int_{0}^{\infty}\frac{r^{k''+1}}{e^{r}}dr \nonumber \\
                                            & = & (k''+1)!=(k-s+1)!=\frac{1}{(k+1)k\cdots (k-s+2)}\cdot (k+1)!. \nonumber  
\end{eqnarray}
Note that
\begin{eqnarray*}
  &      & \left\{ \frac{1}{(k+1)k\cdots (k-s+2)}/\frac{1}{k(k-1)\cdots (k-s+1)}\right\} \\
  & =    & \frac{k(k-1)\cdots (k- s+1)}{(k+1)k\cdots (k-s+2)} \\
  & =    & \frac{(k-s+1)}{k+1} \\ 
  & \leq &1.
\end{eqnarray*}
Let $\rho =1$. Then
$$\int_{0}^{\infty}\frac{r^{k''+1}}{e^{r}}dr\leq \Delta_{0}\cdot \rho^{k+1}\cdot (k+1)!$$
if $\int_{0}^{\infty}\frac{r^{k''}}{e^{r}}dr\leq \Delta_{0}\cdot \rho^{k}\cdot k!$. The smallest $k$ is $k_{0}=2m+2$. Hence the smallest $k''$ is $k_{0}''=k_{0}-s=2m+2-s$. Now,
\begin{eqnarray}
  \int_{0}^{\infty}\frac{r^{k_{0}''}}{e^{r}}dr & = & k_{0}''!=(2m+2-s)!=(k_{0}-s)! \nonumber \\
                                               & = & \frac{1}{k_{0}(k_{0}-1)\cdots (k_{0}-s+1)}\cdot (k_{0}!) \nonumber \\
                                               & = & \Delta_{0}\cdot \rho^{k_{0}}\cdot k_{0}! \ where\ \Delta_{0}=\frac{1}{(2m+2)(2m+1)\cdots (2m-s+3)}. \nonumber 
\end{eqnarray}
It follows that $\int_{0}^{\infty}\frac{r^{k'}}{e^{r}}dr\leq \Delta_{0}\cdot \rho^{k}\cdot k!$ for all $k\geq 2m+2$.

Case3. Assume $k''=k$. Then 
$$\int_{0}^{\infty}\frac{r^{k'}}{e^{r}}dr\leq \int_{0}^{\infty}\frac{r^{k''}}{e^{r}}dr=k!\ \ and\ \ \int_{0}^{\infty}\frac{r^{k'+1}}{e^{r}}dr\leq (k+1)!.$$
Let $\rho=1$. Then $\int_{0}^{\infty}\frac{r^{k'}}{e^{r}}dr\leq \Delta_{0}\cdot \rho^{k}\cdot k!$ for all $k$ where $\Delta_{0}=1$.

The lemma is now an immediate result of the three cases. \hspace{2.7cm} $\sharp$ \\
\\
{\bf Remark}: For the convenience of the reader we should express the constants $\Delta_{0}$ and $\rho$ in a clear form. It's easily shown that\\
(a)$k''>k$ if and only if $n-\lambda>3$,\\
(b)$k''<k$ if and only if $n-\lambda\leq 1$, and\\
(c)$k''=k$ if and only if $1<n-\lambda \leq 3$,\\
where $k''$ and $k$ are as in the proof of the lemma. We thus have the following situations.\\
(a)$n-\lambda>3$. Let $s=\lceil \frac{n-\lambda-3}{2}\rceil $. Then 
$$\rho=1+\frac{s}{2m+3}\ \ and\ \ \Delta_{0}=\frac{(2m+2+s)(2m+1+s)\cdots (2m+3)}{\rho^{2m+2}}.$$
(b)$n-\lambda\leq 1$. Let $s=-\lceil \frac{n-\lambda-3}{2}\rceil $. Then 
$$\rho=1\ \ and\ \ \Delta_{0}=\frac{1}{(2m+2)(2m+1)\cdots (2m-s+3)}.$$
(c)$1<n-\lambda\leq 3$. We have 
$$\rho=1\ \ and \ \ \Delta_{0}=1.$$

Before introducing our main theorem, we need the following two lemmas, first of which is cited directly from\cite{MN3}.
\begin{lem}
  Let $Q,\ Y$, and $\gamma_{n}$ be as in Lemma1.1. Then, given a point $x$ in $Q$, there is a measure $\sigma$ supported on $Y$ such that 
$$\int_{R^{n}}p(y)d\sigma(y)=p(x)$$
for all $p$ in $P_{k}$, and 
$$\int_{R^{n}}d|\sigma|(y)\leq e^{2n\gamma_{n}(k+1)}.$$
\end{lem}
\begin{lem}
  For any positive integer $k$,
$$\frac{\sqrt{(2k)!}}{k!}\leq 2^{k}.$$
\end{lem}
{\bf Proof}. This inequality holds for $k=1$ obviously. We proceed by induction.
\begin{eqnarray*}
  \frac{\sqrt{[2(k+1)]!}}{(k+1)!} &  =   & \frac{\sqrt{(2k+2)!}}{k!(k+1)}=\frac{\sqrt{(2k)!}}{k!}\cdot \frac{\sqrt{(2k+2)(2k+1)}}{k+1}  \\
                                  & \leq & \frac{\sqrt{(2k)!}}{k!}\cdot \frac{\sqrt{(2k+2)^{2}}}{k+1}\leq 2^{k}\cdot \frac{(2k+2)}{k+1}=2^{k+1}. \hspace{2cm}  \ \ \sharp 
\end{eqnarray*}

Because of the local nature of the result, we first restrict our attention to the case where $x$ lies in a cube.
\begin{thm}
  Suppose $h$ is defined as in (4). Let $\mu$ be its corresponding measure as in (6). Then, given a positive number $b_{0}$, there are positive constants $\delta_{0}$ and $\omega$, $0<\omega<1$, which depend on $b_{0}$ for which the following is true:\\
If $f\in {\cal C}_{h,m}$ and $s$ is the $h$ spline that interpolates $f$ on a subset $X$ of $R^{n}$, then
\begin{equation}
  |f(x)-s(x)|\leq \sqrt{l(\lambda,n)}\cdot (\frac{\pi}{2})^{1/4}\cdot \sqrt{n\cdot \alpha_{n}}\cdot c^{\frac{\lambda }{2}}\cdot \sqrt{\Delta_{0}}\cdot \omega^{\frac{1}{\delta}}\cdot \| f\| _{h} \label{eq:main}
\end{equation} 
holds for all $x$ in a cube $E$ provided that (a)$E$ has side $b$ and $b\geq b_{0}$, (b)$0<\delta\leq \delta_{0}$ and (c)every subcube of $E$ of side $\delta$ contains a point of $X$. Here, $l(\lambda,n)$ is defined in (8), $\alpha_{n}$ denotes the volume of the unit ball in $R^{n}$, and $c,\ \Delta_{0}$ are as in (9).

The number $\delta_{0}$ and $\omega$ can be expressed specifically as 
$$\delta_{0}=\frac{1}{3C\gamma_{n}(m+1)},\ \ \omega=\left( \frac{2}{3}\right) ^{\frac{1}{3C\gamma_{n}}}$$
where
$$C=max\left\{ 2\rho'\sqrt{n}e^{2n\gamma_{n}},\ \frac{2}{3b_{0}}\right\} ,\ \ \rho'=\frac{\rho}{c}.$$
The number $\rho$ can be found in the remark following Lemma2.1, $\gamma_{n}$ is defined in Lemma1.1, and $m=1+\frac{\lambda}{2}$ is defined in (4). 
\end{thm}
{\bf Proof}. First, let $\rho,\ \gamma_{n}$, and $b_{0}$ be the constants appearing in Lemma2.1, Lemma1.1, and Theorem2.4, repectively. Let
$$B=2\rho'\sqrt{n}e^{2n\gamma_{n}}\ \ and\ \ C=max\left\{ B,\ \frac{2}{3b_{0}}\right\} $$
where $\rho'=\frac{\rho}{c}$. Let
$$\delta_{0}=\frac{1}{3C\gamma_{n}(m+1)},$$
where $m$ is the order of c.p.d. of $h$.

Now, let $x$ be any point of the cube $E$ and recall that Theorem4.2 of \cite{MN2} implies that
\begin{equation}
  |f(x)-s(x)|\leq c_{k}\| f\| _{h}\int_{R^{n}}|y-x|^{k}d|\sigma|(y)   \label{eq:main1}
\end{equation}
whenever $k>m$, where $\sigma$ is any measure supported on $X$ such that 
\begin{equation}
  \int_{R^{n}}p(y)d\sigma(y)=p(x)  \label{eq:main2}
\end{equation}
for all polynomials $p$ in $P_{k-1}$. Here
$$c_{k}=\left\{ \int_{R^{n}}\frac{|\xi|^{2k}}{(k!)^{2}}d\mu(\xi)\right\} ^{1/2}$$
whenever $k>m$. By (9), for all $2k\geq2m+2$,
\begin{eqnarray}
  c_{k} & =    & \left\{ \int_{R^{n}}\frac{|\xi|^{2k}}{(k!)^{2}}d\mu(\xi)\right\} ^{1/2} \nonumber \\
        & \leq & \frac{1}{k!}\cdot \sqrt{l(\lambda,n)}\cdot (\frac{\pi}{2})^{1/4}\cdot \sqrt{n\alpha_{n}}\cdot c^{-k+\frac{\lambda}{2}}\cdot \sqrt{\Delta_{0}}\cdot \rho^{k}\cdot \sqrt{(2k)!} \nonumber \\
        & \leq & \sqrt{l(\lambda,n)}\cdot (\frac{\pi}{2})^{1/4}\cdot \sqrt{n\alpha_{n}}\cdot c^{\frac{\lambda}{2}}\cdot c^{-k}\cdot \sqrt{\Delta_{0}}\cdot (2\rho)^{k}  
\end{eqnarray}
due to Lemma2.3.

To obtain the desired bound on $|f(x)-s(x)|$, it suffices to find a suitable bound for 
$$I=c_{k}\int_{R^{n}}|y-x|^{k}d|\sigma|(y).$$
This is done by choosing the measure $\sigma$ appropriately. We proceed as follows:

Let $\delta$ be a parameter as in the statement of the theorem. Since $\delta\leq \delta_{0}$ and $0<3C\gamma_{n}\delta\leq \frac{1}{m+1}$, we may choose an integer $k\geq m+1$ so that 
$$1\leq 3C\gamma_{n}k\delta\leq 2.$$
Note that $\gamma_{n}k\delta \leq b_{0}$ for such a $k$. Let $Q$ be any cube which contains $x$, has side $\gamma_{n}k\delta$, and is contained in $E$. Subdivide $Q$ into $(\gamma_{n}k)^{n}$ congruent subcubes of side $\delta$. Since each of these subcubes must contain a point of $X$, select a point of $X$ from each such subcube and call the resulting discrete set $Y$. By virtue of Lemma2.2 we may conclude that there is a measure $\sigma$ supported on $Y$ which satisfies (12) and enjoys the estimate 
\begin{equation}
  \int_{R^{n}}d|\sigma|(y)\leq e^{2n\gamma_{n}k}.  \label{eq:main4}
\end{equation}
We use this measure in (11) to obtain an estimate on $I$.

Using (13), (14), and the fact that support of $\sigma$ is contained in $Q$ whose diameter is $\sqrt{n}\gamma_{n}k\delta$, we may write 
\begin{eqnarray}
  I & \leq & \sqrt{l(\lambda,n)}\cdot (\frac{\pi}{2})^{1/4}\cdot \sqrt{n\alpha_{n}}\cdot c^{\frac{\lambda}{2}}\cdot c^{-              k}\cdot \sqrt{\Delta_{0}}\cdot (2\rho)^{k}\cdot (\sqrt{n}\gamma_{n}k\delta)^{k}e^{2n\gamma_{n}k} \nonumber \\
    & \leq & (C\gamma_{n}k\delta)^{k}(\sqrt{l(\lambda,n)}\cdot (\frac{\pi}{2})^{1/4}\cdot \sqrt{n\alpha_{n}}\cdot c^{\frac{\lambda}{2}}\cdot \sqrt{\Delta_{0}}).
\end{eqnarray}
Since
$$C\gamma_{n}k\delta\leq \frac{2}{3}\ \ and \ \ k\geq \frac{1}{3C\gamma_{n}\delta},$$
(15) implies that 
$$I\leq  \left[ \left( \frac{2}{3}\right) ^{\frac{1}{3C\gamma_{n}}}\right] ^{\frac{1}{\delta}}\cdot \left( \sqrt{l(\lambda,n)}\cdot (\frac{\pi}{2})^{1/4}\cdot \sqrt{n\alpha_{n}}\cdot c^{\frac{\lambda}{2}}\cdot \sqrt{\Delta_{0}}\right) .$$
Hence we may conclude that
$$|f(x)-s(x)|\leq \sqrt{l(\lambda,n)}\cdot (\frac{\pi}{2})^{1/4}\cdot \sqrt{n\alpha_{n}}\cdot c^{\frac{\lambda}{2}}\cdot \sqrt{\Delta_{0}}\cdot \omega^{\frac{1}{\delta}}\cdot \| f\| _{h}$$
where $$\omega=\left( \frac{2}{3}\right) ^{\frac{1}{3C\gamma_{n}}}.$$
This completes the proof. \hspace{8cm} \ \ \ \ $\sharp$

What's noteworthy is that in Theorem2.4 the parameter $\delta$ is not the generally used fill-distance. For easy use we should transform the theorem into a statement described by the fill-distance.

Let$$d(\Omega, X)=sup_{y\in \Omega}inf _{x\in X}|y-x|$$
be the fill-distance. Observe that every cube of side $\delta$ contains a ball of radius $\frac{\delta}{2}$. Thus the subcube condition in Theorem2.4 is satisfied when $\delta=2d(E,X)$. More generally, we can easily conclude the following:
\begin{cor}
  Suppose $h$ is defined as in (4). Let $\mu$ be its corresponding measure as in (6). Then, given a positive number $b_{0}$, there are positive constants $d_{0}$ and $\omega'$, $0<\omega'<1$, which depend on $b_{0}$ for which the following is true: If $f\in {\cal C}_{h,m}$ and $s$ is the $h$ spline that interpolates $f$ on a subset $X$ of $R^{n}$, then
\begin{equation}
  |f(x)-s(x)|\leq \sqrt{l(\lambda,n)}\cdot (\frac{\pi}{2})^{\frac{1}{4}}\cdot \sqrt{n\alpha_{n}}\cdot c^{\frac{\lambda}{2}}\cdot \sqrt{\Delta_{0}}\cdot (\omega')^{\frac{1}{d}}\cdot \| f\| _{h} \label{eq:last}
\end{equation}
holds for all $x$ in a cube $E\subseteq \Omega$, where $\Omega$ is a set which can be expressed as the union of rotations and translations of a fixed cube of side $b_{0}$, provided that (a)$E$ has side $b\geq b_{0}$, (b)$0<d\leq d_{0}$ and (c)every subcube of $E$ of side $2d$ contains a point of $X$. Here, $\alpha_{n}$ denotes the volume of the unit ball in $R^{n}$ and $c,\ \Delta_{0}$ are as in (9). Moreover $d_{0}=\frac{\delta_{0}}{2}$ and $\omega'=\sqrt{\omega}$ where $\delta_{0}$ and $\omega$ are as in Theorem2.4.
\end{cor}
{\bf Proof}. Let $d_{0}=\frac{\delta_{0}}{2}$ and $\delta=2d$. Then $0<d\leq d_{0}$ iff $0<\delta\leq \delta_{0}$. Our corollary follows immediately by noting that $\omega^{\frac{1}{\delta}}=\omega^{\frac{1}{2d}}=\sqrt{\omega}^{\frac{1}{d}}=(\omega')^{\frac{1}{d}}$. \ \ \ $\sharp $\\
\\
{\bf Remark}: The space ${\cal C}_{h,m}$ probably is unfamiliar to most people. It's introduced by Madych and Nelson in \cite{MN1} and \cite{MN2}. Later Luh made characterizations for it in \cite{Lu1} and \cite{Lu2}. Some people think that it's defined by Gelfand and Shilov's generalized Fourier transform, and is therefore difficult to deal with. This is not true. In fact, it can be characterized by Schwartz's generalized Fourier transform. The situation is not so bad. Moreover, some people think that ${\cal C}_{h,m}$ is the closure of Wu and Schaback's function space which is defined in \cite{WS}. This is also not true. The two spaces have very little connection. Luh has also made a clarification for this problem. For further details, please see \cite{Lu3} and \cite{Lu4}.

\end{document}